\documentclass{article}

\usepackage{amsmath,amsthm,amssymb}

\newcommand{\proofend}{\hfill $\blacksquare$}

\newtheorem{theorem}{Theorem}[section]

\newtheorem{corollary}[theorem]{Corollary}

\newtheorem{conjecture}[theorem]{Conjecture}

\input epsf

\title{A note on a conjecture of Gy\'arf\'as} \author{Ryan
R. Martin\footnote{Department of Mathematics, Iowa State
University, Ames, IA 50011.  The author partially supported by the
Clay Mathematics Institute. email: {\tt
rymartin@math.iastate.edu}}}\date{}

\begin{document}

\maketitle

\begin{abstract}
This note proves that, given one member, $T$, of a particular family
of radius-three trees, every radius-two, triangle-free graph, $G$, with
large enough chromatic number contains an induced copy of $T$.
\end{abstract}

\section{Introduction}
A ground-breaking theorem by Erd\H{o}s~\cite{E1} states that for any
positive integers $\chi$ and $g$, there exists a graph with chromatic
number at least $\chi$ and girth at least $g$.  This has an important
corollary.  Let $H$ be a fixed graph which contains a cycle and let
$\chi_0$ be a fixed positive integer.  Then there exists a $G$ such
that $\chi(G)>\chi_0$ and $G$ does not contain $H$ as a subgraph. \\

Gy\'arf\'as~\cite{G1} and Sumner~\cite{S1} independently conjectured
the following:
\begin{conjecture}
For every integer $k$ and tree $T$ there is an integer $f(k,T)$ such
that every $G$ with
\[ \omega(G)\leq k \qquad\mbox{ and }\qquad \chi(G)\geq f(k,T) \]
contains an induced copy of $T$. \\
\label{conKIERSTEAD}
\end{conjecture}

Of course, an acyclic graph need not be a tree.  But,
Conjecture~\ref{conKIERSTEAD} is the same if we replace $T$, by $F$
where $F$ is a forest.  Suppose $F=T_1+\cdots+T_p$ where each $T_i$ is
a tree, then we can see by induction on both $k$ and $p$ that
\[ f(k,F)\leq 2p+|V(F)|f(k-1,F)+\max_{1\leq i\leq
   p}\left\{f(k,T_i)\right\} . \]
A similar proof is given in~\cite{KP1}.  Thus, it is sufficient to
prove Conjecture~\ref{conKIERSTEAD} for trees, as stated. \\

\subsection{Current Progress}
\begin{figure}
\centerline{\epsfbox{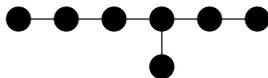}} \caption{Kierstead-Penrice's
$T$} \label{kier-figure1}
\end{figure}

The first major progress on this problem came from Gy\'arf\'as,
Szemer\'edi and Tuza~\cite{GST1} who proved the case when $k=3$ and
$T$ is either a radius two tree or a so-called ``mop.''  A mop is a
graph which is path with a star at the end.  Kierstead and
Penrice~\cite{KP1}  proved the conjecture for $k=3$ and when $T$ is
the graph in Figure~\ref{kier-figure1}. \\

\begin{figure}
\centerline{\epsfbox{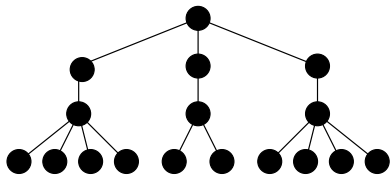}} \caption{A radius three tree
covered in~\cite{KZ1}.} \label{kier-figure2}
\end{figure}

The breakthrough for $k>3$ came through Kierstead and
Penrice~\cite{KP2}, where they proved that
Conjecture~\ref{conKIERSTEAD} is true if $T$ is a radius two tree and
$k$ is any positive integer.  This result contains the one
in~\cite{GST1}.  Furthermore, Kierstead and
Zhu~\cite{KZ1} prove the conjecture true for a certain class of radius
three trees.  These trees are those with all vertices adjacent to the
root having degree 2 or less.  A good example of such a tree is in
Figure~\ref{kier-figure2}.
The paper~\cite{KZ1} contains the result in~\cite{KP1}. \\

Scott~\cite{S2} proved the following theorem:
\begin{theorem}[Scott]
For every integer $k$ and tree $T$ there is an integer $f(k,T)$ such
that every $G$ with $\omega(G)\leq k$ and $\chi(G)\geq f(k,T)$
contains a subdivision of $T$ as an induced subgraph.
\label{tSCOTT}
\end{theorem}
Theorem~\ref{tSCOTT} results in an easy corollary:
\begin{corollary}[Scott]
Conjecture~\ref{conKIERSTEAD} is true if $T$ is a subdivision of a
star and $k$ is any positive integer. \\
\label{cSCOTT}
\end{corollary}

Kierstead and Rodl~\cite{KR1} discuss why
Conjecture~\ref{conKIERSTEAD} does not generalize well to directed
graphs.

\section{The Theorem}
\begin{figure}
\centerline{\epsfbox{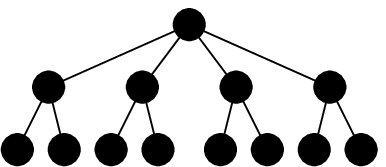}} \caption{$T(4,2)$}
\label{kier-figure3}
\end{figure}

\begin{figure}
\centerline{\epsfbox{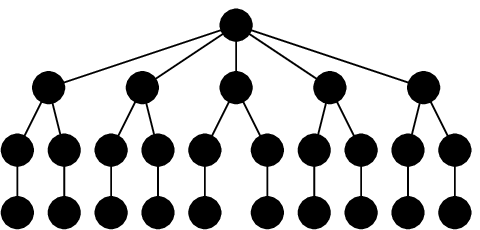}} \caption{$T(5,2,1)$}
\label{kier-figure4}
\end{figure}

In order to prove the theorem, we must define some specific trees.
In general, let $T(a,b)$ denote the radius two tree in which the
root has $a$ children and each of those children itself has
exactly $b$ children.  (Thus, $T(a,b)$ has $1+a+ab$ vertices.) In
particular, $T(t,2)$ is the radius two tree for which the root has
$t$ children and each neighbor of the root has 2 children.
Figure~\ref{kier-figure3} gives a drawing of $T(4,2)$.
Let $T(t,2,1)$ be the radius three tree in which the root has $t$
children, each neighbor of the root has 2 children, each vertex at
distance two from the root has 1 child and each vertex at distance
three from the root
is a leaf.  Figure~\ref{kier-figure4} gives a drawing of $T(5,2,1)$. \\

This allows us to state the theorem:
\begin{theorem}
Let $t$ be a positive integer.  There exists a function $f$, such that
if $G$ is a radius two graph with no triangles and $\chi(G)>f(t)$,
then $G$ must have $T(t,2,1)$ as an induced subgraph.
\label{tKIERSTEAD}
\end{theorem}

\noindent{\bf Proof.}  We will let $r$ be the root of $G$ and let
$S_1=S(r,1)$ be the neighbors of $r$ and $S_2=S(r,2)$ be the
second neighborhood of $r$.  We will try to create a $T(t,2,1)$
with a root $r$ vertex by vertex.  We look for a $v_1\in S_1$ with
the property that there exist $w_{1a},w_{1a}\in N_{S_2}(v_1)$ as
well as $x_{1a}\in N_{S_2}(w_{1a})\setminus
N_{S_2}(w_{1b})\neq\emptyset$ and $x_{1b}\in
N_{S_2}(w_{1b})\setminus N_{S_2}(w_{1a})\neq\emptyset$ such that
$x_{1a}\not\sim x_{1b}$. So, clearly,
$\{v_1,w_{1a},w_{1b},x_{1a},x_{1b}\}$ induce the tree $T(2,1)$.
Let us remove the following vertices from $G$ to create $G_2$:
\[ \{v_1,w_{1a},w_{1b},x_{1a},x_{1b}\}\cup N_{S_2}(v_1)\cup
   N(w_{1a})\cup N(w_{1b})\cup N(x_{1a})\cup N(x_{1b}) . \]
Since $G$ has no triangles, the graph induced by these vertices
has chromatic number at most 4.\footnote{One such coloring is (1)
$N_{S_2}(w_{1a})\cup N_{S_2}(x_{1a})$, (2) $N_{S_2}(w_{1b})\cup
N_{S_2}(x_{1b})$, (3) $N_{S_2}(v_1)$ and (4)
$\{v_1,x_{1a},x_{1b}\}$.} Thus,
$\chi(G_2)\geq\chi(G)-4$. \\

We continue to find $v_2,\ldots,v_s$ from each of $G_2,\ldots,G_s$
in the same manner with $s<t$ so that $G$ has an induced
$T(s,2,1)$ rooted at $r$.  We also have a $G_{s+1}$ so that
$\chi(G_{s+1})\geq \chi(G)-4s$.  If we can continue this process
to the point that $s=t$, we have our $T(t,2,1)$ rooted at $r$. So,
let us suppose that the process stops for some $s<t$. From this
point forward, $S_1$ will actually denote $S_1\cap V(G_{s+1})$ and
$S_2$
will denote $S_2\cap V(G_{s+1})$. \\

Furthermore, in the graph $G_{s+1}$, each vertex $v_1\in S_1$ has
the following property: For any $w_{1a},w_{1b}\in N(v_1)$, the
pair
$$ \left(N_{S_2}(w_{1a})\setminus N_{S_2}(w_{1b}),
         N_{S_2}(w_{1b})\setminus N_{S_2}(w_{1a})\right) $$
induces a complete bipartite graph.  If this were not the case,
then we could find the $x_{1a}$ and $x_{1b}$ that we need. \\

Consider this property in reverse.  Let $v\in S_1$ and $z_1,z_2\in
S_2\setminus N_{S_2}(v)$.  Then the two sets $N_{S_2}(v)\cap
N(z_1)$ and $N_{S_2}(v)\cap N(z_2)$ have the property that one is
inside the other or they are disjoint.  As a result, $N_{S_2}(v)$
has two nonempty subsets such that any $z\in S_2\setminus
N_{S_2}(v)$ has the property that $N_{S_2}(v)\cap N(z)$ contains
either one subset or the other. \\

So, for each $v\in S_2$, there exists some (not necessarily unique
and not necessarily distinct) pair of vertices, $w_a(v),w_b(v)\in
N_{S_2}(v)$ such that for all $z\in S_2$, if $z$ is adjacent to
some member of $N_{S_2}(v)$ then either $z\sim w_a(v)$ or $z\sim
w_b(v)$ or both. \\

For every $v\in S_1$, find such vertices and label them,
arbitrarily as $w_a(v)$ or $w_b(v)$, recognizing that a vertex can
have many labels. Now form the graph $H^{*}$ induced by vertices
from among those labelled as some $w_a(v)$ or $w_b(v)$. Find a
minimal induced subgraph $H$ so that if $h^{*}\in V(H^{*})$, then
there exists $h\in V(H)$ such that $N_{S_2}(h^{*})\subseteq
N_{S_2}(h)$.
\\

We have a series of claims that end the proof: \\

\noindent{\bf Claim 1.} $\chi(H)=\chi(S_2)$. \\
\noindent{\bf Proof of Claim 1.} Since $H$ is a subgraph of $S_2$,
$\chi(H)\leq\chi(S_2)$. If we properly color $H$ with $\chi(H)$
colors, then we can extend this to a coloring of $S_2$.  We do
this by giving $z\in S_2$ the same color as that of some $h\in
V(H)$ with the property that $N_{S_2}(z)\subseteq N_{S_2}(h)$.

This is possible first because there must be some $h^{*}=w_A(v)$
or $h^{*}=w_B(v)$ in $H^{*}$ with $N_{S_2}(z)\subseteq
N_{S_2}(h^{*})$.  Further, there is an $h$ such that
$N_{S_2}(h^{*})\subseteq N_{S_2}(h)$.  So, $N_{S_2}(z)\subseteq
N_{S_2}(h)$.  Now suppose $z_1$ and $z_2$ are given the same color
but are adjacent.  Let $h_1$ and $h_2$ be the vertices in $H$
whose neighborhoods dominate those of $z_1$ and $z_2$,
respectively and whose colors $z_1$ and $z_2$ inherit.  Because
$z_1\sim z_2$, $h_1\sim z_2$ and $h_2\sim z_1$.  But then it must
also be the case that $h_1\sim h_2$.  Thus, $h_1$ and $h_2$ cannot
receive the same color, a contradiction. \proofend \\

\noindent{\bf Claim 2.} $H$ induces a $T(2t+1,8)$. \\
\noindent{\bf Proof of Claim 2.} Because $S_1$ is an independent
set, $\chi(S_2)\geq\chi(G_{s+1})-1$.  Because $\chi(G)$, hence
$\chi(G_{s+1})$, is large, Claim 1 ensures that $\chi(H)$ is
large. Claim 2 results from \cite{GST1}, because $T(2t+1,8)$ is a
radius-two tree.
\proofend \\

\begin{figure}
\centerline{\epsfbox{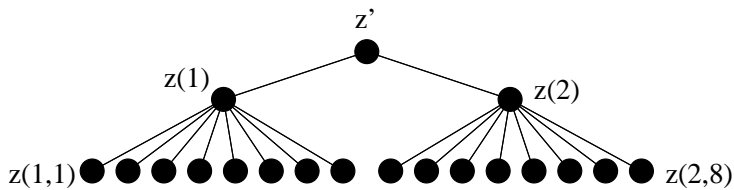}} \caption{$T(2,8)$ with some
vertices labelled} \label{kier-figure5}
\end{figure}

Let the tree $T$, guaranteed by Claim 2, have root $z'$, its
children be labelled $z(1),\ldots,z(2t+1)$ and the children of
each $z(i)$ be labelled $z(i,1),\ldots,z(i,8)$.  Figure
\ref{kier-figure5} shows one such tree. \\

\noindent{\bf Claim 3.} If $v\in S_1$ is adjacent to $z(i,j)$,
then $v$ cannot be adjacent to any other vertices of $T$ except
one other vertex $z(i,j')$ or $z'$. \\
\noindent{\bf Proof of Claim 3.}  If $v\in S_1$ is adjacent to,
say, $z(1,1)$, then $v\not\sim z(i,j)$ if $i\neq 1$.  This is
because $N_{S_2}(w_A(v))\triangle N_{S_2}(w_B(v))$ induces a
complete bipartite graph which would imply an edge between $z(1)$
and $z(i)$.

It can be shown, for similar reasons, that if $v\sim z(1,1)$, then
$v\not\sim z(i)$ for any $i\neq 1$.  Also, $v\not\sim z(1)$
because $G$ is triangle-free.  \proofend \\

\noindent{\bf Claim 4.} We may assume that there is a $v_1\in S_1$
that is adjacent to (without loss of generality) $z(1,1)$ as well
as $z'$. \\
\noindent{\bf Proof of Claim 4.} We prove this by contradiction.
Applying Claim 3 to every leaf of $T$, we see that since Claim 4
is not true, then for $i=1,\ldots,2t+1$, we can find a set of $4$
vertices of the form $z(i,j)$ and $4$ vertices from $S_1$ so that
they induce a perfect matching.  Furthermore, the $4(2t+1)$
vertices from $S_1$ are each adjacent to no other vertices of $T$,
because of Claim 3. Hence, we have our induced $T(t,2,1)$, a
contradiction.
\proofend \\

Because our definition of $H$ guaranteed that vertices had
neighborhoods that were not nested, there must be some $z''\in
S_2$ that is adjacent to $z(1,1)$ but not $z'$. Call this vertex
$z''$. \\

\noindent{\bf Claim 5.} For any $z(i,j)$ with $i\neq 1$ and any
$v\in S_1$ adjacent to $z(i,j)$, $v$ cannot be adjacent to both
$z'$ and $z''$. \\
\noindent{\bf Proof of Claim 5.} We again proceed by
contradiction, supposing that $v\sim z(i,j),z',z''$.  There is,
without loss of generality, $w_a(v)\in N_{S_2}(v)$ such that
$N_{S_2}(z'')\subseteq N_{S_2}\left(w_a(v)\right)$. Thus, either
$N_{S_2}(z')\subseteq N_{S_2}\left(w_a(v)\right)$ or
$N_{S_2}\left(z(i,j)\right)\subseteq N_{S_2}\left(w_a(v)\right)$.
But if $w_a(v)$ were deleted from $H^{*}$ to form $H$, either $z'$
or $z(i,j)$ would have been deleted as well.

Therefore, either $w_a(v)=z'$ or $w_a(v)=z(i,j)$.  So,
$N_{S_2}(z'')\subseteq N_{S_2}(z')$ or $N_{S_2}(z'')\subseteq
N_{S_2}\left(z(i,j)\right)$.  We can conclude that either $z'\sim
z(1,1)$ or $z(i,j)\sim z(1,1)$.  This contradicts the fact that
$T$ is an induced subtree. \proofend \\

\noindent{\bf Claim 6.} For all $i\neq 1$, $z''$ is adjacent to
$z(i)$ but no vertex $z(i,j)$. \\
\noindent{\bf Proof of Claim 6.} Note that $z(2),\ldots,z(2t+1)$
are adjacent to $z'$ but not $z(1,1)$.  Because of the condition
that $N_{S_2}(z')\triangle N_{S_2}\left(z(1,1)\right)$ induces a
complete bipartite graph, $z''$ must be adjacent to
$z(2),\ldots,z(2t+1)$.  Because $G$ is triangle-free, $z''$ cannot
be adjacent to any vertex of the form $z(i,j)$ where $i\neq 1$.
\proofend \\

Now we construct the tree we need.  For each $z(i,j)$, $i\neq 1$,
find a vertex $v(i,j)\in S_1$ to which $z(i,j)$ is adjacent.
According to Claim 3, no $v(i,j)$ vertex can be adjacent to any
vertex of $V(T)\setminus\{z'\}$ and, according to Claim 5, it is
adjacent to at most one of $\{z',z''\}$. \\

For each $i\in\{2,\ldots,2t+1\}$, the majority of
$\{v(i,1),\ldots,v(i,8)\}$ have that $v(i,j)$ is either
nonadjacent to $z'$ or nonadjacent to $z''$.  Without loss of
generality, we conclude that $z'$ has the property that, for
$i=2,\ldots,t+1$, the vertices $v(i,1),\ldots,v(i,4)$ fail to be
adjacent to $z'$. \\

Since any vertex of $S_1$ can be adjacent to at most two vertices
of $H$, then for $i=2,\ldots,t+1$,
$\left|\{v(i,1),\ldots,v(i,4)\}\right|\geq 2$. Therefore, we
assume that for each $i\in\{2,\ldots,t+1\}$, $v(i,1)$ and $v(i,2)$
are distinct.  But now the vertex set
$$ \{z'\}\cup\bigcup_{i=2}^{t+1}
   \left(\left\{z(i),z(i,1),z(i,2),v(i,1),v(i,2)\right\}\right) $$
induces $T(t,2,1)$. \proofend

\section{Acknowledgements}

Very many thanks to an anonymous referee who identified a
significant error in the first version of this paper.

\bibliography{kierstead}
\bibliographystyle{plain}

\end{document}